\documentstyle[12pt]{amsart}

\newcommand{\sab}{\left(\begin{array}{rr} a & b\\-b & a \end{array}\right)}
\newcommand{\ts}{\vspace{\baselineskip}\noindent{\bf Proof.}$\;\;$}
\newcommand{\GG}{{\bf G}}
\newcommand{\ZZ}{{\bf Z}}
\newcommand{\QQ}{{\bf Q}}
\newcommand{\RR}{{\bf R}}
\newcommand{\CC}{{\bf C}}
\newcommand{\PP}{{\bf P}}

\begin{document}

\title[Kuga-Satake varieties and the Hodge conjecture]{Kuga-Satake varieties and the Hodge conjecture}

\author{Bert van Geemen}
\address{Dipartimento di Matematica, Universit\`a di Pavia, 
Via Ferrata 1, I-27100 Pavia, Italia}
\email{geemen@@dragon.ian.pv.cnr.it}

\maketitle
\begin{center}
{\sc Introduction}
\end{center}
\vskip15pt

Kuga-Satake varieties are abelian varieties associated to certain weight two
Hodge structures, for example the second cohomology group of a K3 surface.
We start with an introduction to Hodge structures 
and we give a detailed account of the construction of Kuga-Satake varieties.
The Hodge conjecture is discussed in section \ref{HC}. An excellent survey of the Hodge conjecture for abelian varieties is \cite{BG}.

 We point out a connection between the Hodge conjecture for abelian varieties and Kuga-Satake
varieties in section \ref{Lom}. In section \ref{KSH} we discuss the implications
of the Hodge conjecture on the geometry of surfaces whose 
second cohomology group has a Kuga-Satake variety.
We conclude with some recent results, inspired by an example of C.\ Voisin,
on Kuga-Satake varieties of Hodge structures on which an imaginary quadratic field acts.

I'm indebted to E.\ Izadi, G.\ Lombardo and M.\ Nori for helpful discussions.

\section{Polarized Hodge structures}
\subsection{Definition.} A (rational) Hodge structure of weight $k\;(\in\ZZ)$
is a $\QQ$-vector space $V$ with a decomposition of its complexification $V_\CC:=V\otimes_\QQ\CC$:
$$
V_\CC=\oplus_{p+q=k}V^{p,q},\qquad{\rm and}\quad \overline{V^{p,q}}=V^{q,p}
\qquad(p,\,q\in\ZZ).
$$
Here complex conjugation on $V_\CC$ is given by $\overline{v\otimes z}:=v\otimes\bar{z}$ for $v\in V$ and $z\in\CC$. 
We will take $k,\,p,\,q\geq 0$ except in some of the proofs.

\subsection{}\label{betticoh}
The (Betti) cohomology groups
$H^k(X,\QQ)$ of a complex smooth projective variety are Hodge structures:
$$
H^k(X,\CC)=\oplus_{p+q=k}H^{p,q}(X).
$$
Identifying $H^ k(X,\CC)$ with harmonic differential forms, the subspaces
$H^{p,q}(X)$ consist of the harmonic forms of type $(p,q)$ and
$H^{p,q}(X)\cong H^q(X,\Omega^p)$.

\subsection{}
It is useful to identify Hodge structures on $V$ with certain representations of the group 
$\CC^*$ on $V_\RR:=V\otimes_\QQ\RR$. We identify $\CC^*$ with a subgroup
of $GL(2,\RR)$:
$$
\CC^*\cong\Bigl\{s(a,b):=
\sab\in GL(2,\RR):\;a^2+b^2\neq 0\Bigr\},\qquad
z=a+bi\longmapsto s(a,b).
$$
The eigenvalues of $s(a,b)$ are $z=a+bi$ and $\bar{z}=a-bi$
with corresponding eigenvectors $e_1+ie_2$ and $e_1-ie_2$, here $e_j$ is 
the standard $j$-th basis vector of $\CC^2$.
An algebraic representation of $\CC^*$
is defined to be a homomorphism
$$
h:\CC^*\longrightarrow GL(V_\RR)
$$
such that, with respect to some basis of $V_\RR$, the entries of the matrix
$h(s(a,b))$ are polynomials, with coefficients in $\RR$, in $a,\,b,\,(a^2+b^2)^{-1}$. The following proposition is well-known:

\subsection{Proposition.} There is a bijection between rational Hodge structures
of weight $k$ on a $\QQ$-vector space $V$ and 
algebraic representations $h:\CC^*\rightarrow GL(V_\RR)$ with $h(t)=t^k$ for $t\in\RR$. The Hodge structure defined by $h$ is denoted by $(V,h)$
and its Hodge decomposition is:
$$
V^{p,q}:=\{v\in V_\CC:\;h(z)v=z^p\bar{z}^qv\,\}.
$$

\ts 
Composing $h$ with the inclusion $GL(V_\RR)\subset GL(V_\CC)$
we get a representation of $\CC^*$ on $V_\CC$
and the matrix coefficients of this representation are polynomials in
$z,\,\bar{z}$ and $(z\bar{z})^{-1}$.
There is a basis $\{v_i\}$ of $V_\CC$
of simultaneous eigenvectors: $h(z)v_i=\lambda_i(z)v_i$ for some homomorphisms
$\lambda_i:\CC^*\rightarrow \CC^*$. As $\lambda_i$ is a polynomial in $z,\;\bar{z}$ and $(z\bar{z})^{-1}$ we get $\lambda_i(z)=z^p\bar{z}^q$ for some $p,\,q\in\ZZ$.
Since also the conjugate of an eigenvalue is an eigenvalue (on the conjugated eigenvector) we have a Hodge structure.

Conversely, any element in $V_\CC$ can be written as $v\otimes 1+w\otimes i$ with $v,\,w\in V_\RR$. If $v\otimes 1+w\otimes i\in V^{p,q}$ then $v\otimes 1-w\otimes i\in V^{q,p}$.
Let $\{v_r+iw_r\}_r$ be a basis of $V^{p,q}$ with $p\geq q$ and define $V_p=\langle v_r,\,w_r\rangle_r\;(\subset V_\RR)$ be the span of the $v_r,\,w_r$'s. 
Then 
$$
V_\RR=\oplus_{p\geq q} V_p,\qquad {\rm and}
\quad V_p\otimes_\RR\CC=V^{p,q}\oplus V^{q,p}.
$$
The representation $h$ of $\CC^*$ on $V_{\RR}$ is constructed
on each of the subspaces $V_p$.

In case $p=q$ we have $V^{p,p}=\overline{V^{p,p}}$, 
so $V^{p,p}=V_p\otimes \CC$.
Define $h(a+bi)v:=(a^2+b^2)^pv$ for all $v\in V_p$.

Next fix $p,q$ with $p=q+l$ and $l>0$ and let
$\{v_r +iw_r\}_r$ be a basis for $V^{p,q}$. Then the vectors $\{v_r,\,w_r\}_r\in V_{\RR}$, are independent 
over $\RR$ and $V_p=\oplus_r \langle v_r,\,w_r\rangle$. 
For each $r$ define a representation of $\CC^*$ on 
$\langle v_r,\,w_r\rangle\subseteq V_p$ by 
\[ 
h(a+bi)=(a^2+b^2)^q \cdot \sab^l \; . 
\]
The eigenspaces in $V_\CC$ of this representation 
are the $V^{p,q}$'s. \qed

\subsection{Tensor products of Hodge structures.} \label{tensor}
We transfer the usual algebra constructions on representations to Hodge structures.
Given rational Hodge structures $(V,h_V)$, $(W,h_W)$ of weight $k_V,\,k_W$
one defines a rational Hodge structure $(V\otimes W,h_V\otimes h_W)$
of weight $k_V+k_W$ by:
$$
h_V\otimes h_W:\CC^*\longrightarrow GL((V\otimes W)_\RR),\qquad
z\longmapsto[v\otimes w\mapsto (h_V(z)w)\otimes (h_W(z)w)].
$$
The dual vector space $V^*:=Hom_\QQ(V,\QQ)$ is also a Hodge structure
(of weight $-k_V$):
$$
h^*_V:\CC^*\longrightarrow GL(V_\RR^*),\qquad
(h_V^*(z)f)(v):= f(h_V(z)^{-1}v),
$$
here $f\in V^*_\RR=Hom_\RR(V_\RR,\RR)$ and $v\in V_\RR$.
The Tate Hodge structure $\QQ(n)$ ($n\in \ZZ$) is defined by
the vector space $\QQ$ and the homomorphism:
$$
h_n:\CC^*\longrightarrow GL_1(\RR),\qquad z\longmapsto (z\bar{z})^{-n},
$$
it has weight $-2n$ and $\QQ(n)^{p,q}=0$ unless $p=q=-n$
in which case $\QQ(n)^{-n,-n}=\CC$. We write $V(n):=V\otimes \QQ(n)$,
it is a Hodge structure of weight $k_V-2n$ with $V(n)^{p,q}=V^{p+n,q+n}$.

\subsection{Morphisms of Hodge structures.}
A morphism of Hodge structures $f:(V,h_V)\rightarrow (W,h_W)$ is a linear
map $f:V\rightarrow W$ such that $f$ intertwines the representations
$h_V$ and $h_W$ up to a Tate twist:
$$
f(h_V(z)v)=(z\bar{z})^{n}h_W(z)f(v)
$$ 
($f$ gives a `strict' morphism of Hodge structures $f:V\rightarrow W(-n)$). A morphism of Hodge structures satisfies $f_\CC(V^{p,q})\subset W^{p+n,q+n}$, here $f_\CC$ is the $\CC$-linear extension of $f$.
Since $f$ commutes with the $\CC^*$-representations, both kernel and image of
$f$ are (sub) Hodge structures.

 We denote 
by $Hom_{Hod}(V,W)\;(\subset Hom_\QQ(V,W))$
the $\QQ$-vector space of morphisms of Hodge structures.

\subsection{Definition.} 
Let $V$ be a rational Hodge structure of weight $k$
and let $h:\CC^*\rightarrow GL(V_\RR)$ be the corresponding representation.
A polarization on $V$ is a bilinear map:
\[ 
\Psi: V \times V \longrightarrow \QQ  
\]
satisfying (for all $v,w\in V_{\RR}$): 
\[
\Psi(h(z)v,h(z)w)=(z\bar{z})^k\Psi(v,w)
\]
and 
\[ 
\Psi (v,h(i)w)\;\;\;
\mbox{\rm is a symmetric and positive definite form:}
\] 
$\Psi(v,h(i)w)=\Psi(w,h(i)v)$ for all $v,w\in V_\RR$ and 
$\Psi(v,h(i)v)>0$ for all $v\in V_\RR-\{0\}$.

The map $h(i):V_\RR\rightarrow V_\RR$ is called the Weil operator.
A polarization is a strict morphism of Hodge structures
$V\otimes V\rightarrow \QQ(-k)$ .

\subsection{Lemma.}\label{lempol}
 Let $(V,\Psi)$ be a rational polarized Hodge structure of weight $k$. We denote the $\CC$-linear extension of $\Psi$ by $\Psi_\CC$. Then:
\begin{enumerate}
\item $\Psi$ is symmetric if $k$ is even, and is alternating if $k$ is odd.
\item For $x_{p,q}\in V^{p,q}$ and $y_{r,s}\in V^{r,s}$:
$$
\Psi_\CC(x_{p,q},y_{r,s})=0\qquad {\rm if}\quad (p,q)\neq (s,r).
$$
In particular, the direct sum decomposition:
$$
V_{\RR}=\oplus_{p} V_p,\qquad V_p\otimes_\RR\CC=V^{p,q}\oplus V^{q,p}
$$ 
is an orthogonal direct sum w.r.t.\ $Q$.
\item The (restriction of the) $\CC$-bilinear form 
$$
\Psi_\CC: V^{p,q}\times V^{q,p}\longrightarrow \CC
$$
is non-degenerate (so $V^{p,q}\cong_{\Psi_\CC}\, (V^{q,p})^{dual}$).
\item If the weight is even, the quadratic form
  defined by the $\RR$-linear extension of $\Psi$:
$$
Q:V_\RR\longrightarrow \RR,\qquad Q(v):=\Psi_\RR(v,v),
$$
satisfies:
$
(-1)^{l-p}Q_{|V_p}>0
$ where $k=2l$.
\end{enumerate}

\ts
Note that $h(i)^2 v=h(-1)v=(-1)^k v$
for $v\in V_{\RR}$. Now use:
\[ 
\Psi(v,w)=\Psi((h(i)v),h(i)w)=\Psi(w,h(i)^2v)=(-1)^k\Psi(w,v).
\]
Next we observe that for all $z\in\CC^*$:
$$
(z\bar{z})^k\Psi_\CC(x_{p,q},y_{r,s})=\Psi_\CC(h(z)x_{p,q},h(z)y_{r,s})=
\Psi_\CC(z^p\bar{z}^qx_{p,q},z^r\bar{z}^sy_{r,s})=z^{p+r}\bar{z}^{q+s}
\Psi_\CC(x_{p,q},y_{r,s}).
$$
Hence $\Psi_\CC(x_{p,q},y_{r,s})$ can be non-trivial only when $p+r=q+s=k$.
Since also $p+q=r+s=k$ the second statement follows.

Let $x\in V^{p,p}$, then $x=\lambda v$ with $v\in V_p$, $\lambda\in\CC$
and $h(i)v=v$. 
Then $\Psi_\CC(x,x)=\lambda^2\Psi(v,v)=\lambda^2\Psi(v,h(i)v)$ and
$\Psi(v,h(i)v)$ is non-zero if $v\neq 0$.
If $p\neq q$ then for any non-zero $x\in V^{p,q}$, $\bar{x}\in V^{q,p}$ 
and thus
$x+\bar{x}\in V_p$ is non-zero. Therefore:
$$
\begin{array}{rcl}
0&<& \Psi(x+\bar{x},h(i)(x+\bar{x}))\\
&=&\Psi(x,h(i)\bar{x})+\Psi(\bar{x},h(i)x)
      +0+0\\
&=&2i^{q-p}\Psi(x,\bar{x}),
\end{array}
$$
in the last step we used the symmetry of $\Psi(x,h(i)y)$. 
Thus $\Psi$ gives a non-degenerate pairing.

Since $p+q=k=2l$, $h(i)=i^{p-q}=i^{2p-2l}=(-1)^{p-l}$ on $V^{p,q}$,
and also on $V^{q,p}$, hence $h(i)$ acts as $(-1)^{l-p}$ on $V_p$.
As $\Psi(v,h(i)v)$ is positive definite, 
$Q$ is positive definite on $V_p$ if $l-p$ is even and negative definite otherwise.
\qed

\subsection{Example.} 
The cohomology groups $H^k(X,\QQ)$ (as in \ref{betticoh})
have a polarization
see \cite{WeilP}, Th\'eor\'eme IV.7 and corollaire or \cite{GH}, p.\ 123.

\section{The Hodge conjecture}\label{HC}

\subsection{Hodge cycles.}
The space of Hodge classes in a rational Hodge structure $V$ 
of weight $k$ is the 
$\QQ$-subvector space:
$$
B(V):=\left\{ \begin{array}{ccl}
0&{\rm if}\quad &k\;\;\mbox{is odd},\\
V \cap V^{p,p}\quad&
{\rm if}\quad &k=2p.
\end{array}\right.
$$
Note that $V\hookrightarrow V_\CC=V\otimes_\QQ\CC,\;v\mapsto v\otimes 1$
and the intersection $V \cap V^{p,p}$ takes place in $V_\CC$.

\subsection{Algebraic cycles.} Let $X$ be smooth algebraic variety and let $Z$
be (any) irreducible subvariety of codimension $p$ in $X$. Then $Z$ defines a cohomology class $[Z]\in H^{2p}(X,\QQ)$. This defines a cycle class map from
the Chow group (with coefficients in $\QQ$) of codimension $p$ cycles
$$
[.]:CH^p(X)_\QQ\longrightarrow H^{2p}(X,\QQ),\qquad
\sum a_iZ_i\longmapsto \sum a_i[Z_i].
$$
It is well-known that the image of this map lies in the space of Hodge cycles:
$$
[CH^p(X)_\QQ]\subset B(H^{2p}(X,\QQ))\;=H^{2p}(X,\QQ)\cap H^{p,p}(X).
$$

\subsection{The Hodge conjecture.} The Hodge conjecture asserts that:
$$
[CH^p(X)_\QQ]=B(H^{2p}(X,\QQ))
$$
that is, any cohomology class in $H^{2p}(X,\QQ)$ which is of type $(p,p)$ 
is the class of a codimension $p$ cycle. 

\subsection{}
The Hodge conjecture is known to be true in case $p=1$ (it follows from an analysis of the exponential sequence) and thus is also true if $p=-1+\dim X$
and, obviously also in the cases $p=0,\,\dim X$.
In other cases however the conjecture is very much open.

\subsection{Morphisms of Hodge structures and Hodge classes.}
Let $f:V\rightarrow W$ be a (strict) morphism of Hodge structures. Then 
$f\in Hom_\QQ(V,W)=V^*\otimes W$, which is a Hodge structure. Moreover,
we have, for all $z\in\CC^*$ and all $v\in V_\RR$:
$$
f(h_V(z)v)=h_W(z)f(v)\qquad{\rm so}\quad
f(v)=h_W(z)f(h_V(z)^{-1}v)=((h_V^*(z)\otimes h_W(z))f)(v),
$$
thus $f$ is an invariant in the representation $h_V^*\otimes h_W$ of
$\CC^*$ on $(V^*\otimes W)_\RR$, and hence $f$ is of type $(0,0)$. (In case one has $f(h_V(z)v)=(z\bar{z})^{n}h_W(z)f(v)$, $f$ is of type $(n,n)$.)
Thus a morphism of Hodge structures is a Hodge class
and one finds that:
$$
Hom_{Hod}(V,W)\cong B(V^*\otimes W).
$$

\subsection{Morphisms and correspondences.}\label{morcor}
We apply the relation between morphisms of Hodge structures and Hodge classes
to the cohomology of algebraic varieties.
Let $V\;(\subset H^{k}(X,\QQ))$ 
be a sub-Hodge structure, 
the orthogonal projection $\pi_V:H^{k}(X,\QQ)\rightarrow V$
(use the polarization on $H^k$)
and the inclusion $i_V:V\hookrightarrow H^k(X,\QQ)$
are morphisms of Hodge structures. Let
$W\;(\subset H^l(Y,\QQ))$ be a sub-Hodge structure isomorphic to $V$,
let $g:V\stackrel{\cong}{\rightarrow} W$.
Then
$$
f:=i_Wg\pi_V:H^k(X,\QQ)\longrightarrow H^l(Y,\QQ)
$$
is a morphism of Hodge structures and hence 
$f\in B(H^k(X,\QQ)^*\otimes H^l(Y,\QQ))$.

The cohomology of $X\times Y$ can be computed by the K\"unneth formula:
$$
H^n(X\times Y,\QQ)=\oplus_{p+q=n} H^p(X,\QQ)\otimes H^q(Y,\QQ).
$$
Poincar\'e duality shows that $H^p(X,\QQ)\cong H^{2d-p}(X,\QQ)^*$ where $d=\dim X$. Thus $f$ defines a Hodge class in
$$
f\in B(H^k(X,\QQ)^*\otimes H^l(Y,\QQ))
=B(H^{2d-k}(X,\QQ)\otimes H^l(Y,\QQ)))\hookrightarrow 
B(H^{2d-k+l}(X\times Y,\QQ)).
$$
According to the Hodge conjecture, there should exist a cycle $Z$ on 
$X\times Y$ with cycle class $[Z]=f$. A cycle on $X\times Y$ is also called
a correspondence between $X$ and $Y$.

Thus the Hodge conjecture implies that any morphism of Hodge structures
has a geometric origin in a cycle $Z$ on the product of the varieties.
This observation is of importance for the theory of (the category of
Hodge) motives in which, roughly
speaking, the objects are sub-Hodge structures and the morphisms are
(equivalence classes of) correspondences.

\section{The Mumford-Tate group.}
\subsection{}
The Mumford-Tate group $MT(V)$
of a rational Hodge structure $(V,h)$ is an algebraic subgroup of $GL(V)$. 
It allows one to find the Hodge cycles 
$B(V^{\otimes m})$, for any $m$, using representation theory of Lie groups.

\subsection{}
As the Hodge structures obtained from smooth projective varieties are
polarized, we restrict ourselves to that case.
A polarization $\Psi$ on a weight $k$ Hodge structure $(V,h)$
may be considered as an element of $B(V^*\otimes V^*)$ and the isomorphism
$$
V\stackrel{\cong}{\longrightarrow} V^*,\qquad v\longmapsto[w\mapsto\Psi(v,w)]
$$
is a morphism of Hodge structures
(it is a strict morphism $V\rightarrow V^*(-k)$). 

An algebraic group $G\;(\subset GL(V))$ defined over $\QQ$
is defined by polynomial equations in the $n^2$ matrix entries $g_{ij}$ of
$g$ and $\det(g)^{-1}$
with
coefficients in $\QQ$.
For a commutative $\QQ$-algebra $A$ (for example $A=\QQ,\,\RR,\,\CC$),
the set of solutions to these equations in $A^{n^2}$ is a group
(under matrix multiplication) and this group is denoted by $G(A)$.

The algebraic subgroup of $GL(V)$ `fixing' $\Psi$ is defined to be:
$$
G(\Psi):=\{ g\in GL(V):\; \Psi(gv,gw)=\nu(g)\Psi(v,w)\quad
\mbox{for some}\;\nu(g)\in \GG_m\}.
$$
Here $\GG_m$ is the multiplicative group, so $\GG_m(\QQ)=\QQ^*$ etc.
In case the weight is odd we usually write $GSp(\Psi)$
and if the weight is even
$GO(\Psi)$ for $G(\Psi)$.

\subsection{Definition.} Let $(V,\Psi)$ be a polarized rational Hodge 
structure of weight $k$.
\begin{enumerate}
\item $G_1$ is the algebraic subgroup  of the $g\in G(\Psi)$
  for which there is a $\omega(g)\in\GG_m$ such that
  $  g\cdot t=\omega(g)^{m} t$ for all  $t\in B(V^{\otimes m})$.

\item  $G_2$ is the smallest algebraic subgroup of $GL(V)$ which
is defined over $\QQ$ and which still satisfies 
\[
 h(\CC^*) \subseteq G_2(\RR). 
\]
\end{enumerate}

\subsection{Remark.}
The condition is that $g\in G_1$ acts as a scalar multiple of the identity on
the $\QQ$-subspace $B(V^{\otimes m})$ of $V^{\otimes m}$, the shape of the
scalar follows from tensoring Hodge cycles.
Since $h(\CC^*)\in G(\Psi)$ we have $G_2\subset G(\Psi)$.

\subsection{Theorem.}\label{mt}
For any polarized rational Hodge structure $(V,h,\Psi)$ one has:
\[ 
G_1=G_2=:\, MT(V) \qquad{\rm and}\qquad MT(V)\subset G(\Psi), 
\]
the algebraic group $MT(V)$ is called the Mumford-Tate group of $(V,h,\Psi)$.
The Hodge cycles in the Hodge structure $(V^{\otimes m},h^{\otimes m})$ are
the $MT(V)$-`invariants':
$$
B(V^{\otimes m})=
\{w\in V^{\otimes m}:\; gw=\omega(g)^m w\quad\forall g\in MT(V) \}.
$$

The Mumford-Tate group is reductive
(i.e. every finite dimensional representation of $MT(V)$
is a direct sum of irreducible representations).

\ts
This is all in \cite{DMOS}, but we defined $MT(V)$ as the projection in
$GL(V)$ of the group $G$ in I.3
and we use the polarization to identify $V$ and $V^*$.
Proposition I.3.4 shows that $G_1=G_2$ and 
Proposition I.3.1 identifies the Hodge cycles as $MT(V)$-invariants.
The reductivity is proven with Weyl's unitary trick
in  Proposition I.3.6, again the polarization is essential.
\qed

\subsection{Corollary.}
We have: 
$
End_{Hod}(V)\cong End_{MT}(V)$,
{with}
$$
End_{MT}(V)=\{M\in End(V):\; Mg=gM\quad\forall g\in MT(V)\}.
$$

\ts
This follows from Theorem \ref{mt} since $End(V)=V^*\otimes V
\cong V\otimes V$,
hence $End_{Hod}(V)\cong B(V^{\otimes 2})$ and
$End_{MT}(V)$ corresponds to the space of $MT(V)$-invariants in $V^{\otimes 2}$.
\qed

\section{Hodge structures of weight one and two}

\subsection{}
Given a polarized Hodge structure $(V,h,Q)$ of weight 2, it is
interesting to know
if there exists a Hodge structure $(W,h_W)$ of weigth 1 with
$$
V\hookrightarrow W\otimes W
$$
(we do not even require that $W$ be polarized). 
Note that if there is a pair of weight 1 Hodge structures $W_1,\;W_2$ with
$V\subset W_1\otimes W_2$ then
$V\subset W\otimes W$ with $W=W_1\oplus W_2$.

In general the existence of $W$ depends very much on the Mumford-Tate group 
of $V$ but in case $MT(V)=GO(Q)$, the answer is very easy. 
The proof of the following proposition
follows $\S 7$, Remarques of \cite{DK3}. Deligne's observation has been 
further 
explored by C.\ Schoen in \cite{Scov}. 

\subsection{Proposition.}
Let $(V,h,Q)$ be a polarized Hodge structure of weight 2 with $\dim V^{2,0}>0$
and $MT(V)=GO(Q)$.
  
If 
$V\subset W\otimes W$ for some Hodge structure $W$ of 
weight 1 we must have $\dim V^{2,0}=1$.

\ts
The tangent space at the identity element of an algebraic group is the
Lie algebra of that group. The inclusions $SO(Q)\subset GO(Q)=\GG_m O(Q)
\subset GL(V)$ give inclusion of Lie algebras
$so(Q)\subset go(Q)\subset End(V)$.
The sub-space $so(Q)$ 
 is defined by the linear equations
$\sum_{ij}(\partial{F}/{\partial{g_{ij}}})(I)M_{ij}=0$ where
$F\in \QQ[\ldots,g_{ij},\ldots]$ runs over the equations defining $SO(Q)$.
Since for any $F\in \QQ[\ldots,g_{ij},\ldots]$ we have
$F(\ldots,\delta_{ij}+\epsilon m_{ij},\ldots)=F(I)+\epsilon 
\sum_{ij}(\partial{F}/{\partial{g_{ij}}})(I)M_{ij}$ where $I=(\delta_{ij})$
and $\epsilon^2=0$,
$so(Q)$ is the subspace of 
$M\in End(V)$ such that
$Q((I+\epsilon M)x,(I+\epsilon M)y)=Q(x,y)$, that is:
$$
so(Q)=\{M\in End(V):\; Q(Mx,y)+Q(x,My)=0\quad \forall x,\,y\in V\}.
$$

Since $h(z)\in GO(Q)$ one has $h(z)SO(Q)h^{-1}(z)\subset SO(Q)$,
the differential of this automorphism is called $Ad(h(z))$.
It defines a natural Hodge structure on $so(Q)$: 
$$
Ad(h):\CC^*\longrightarrow GL(so(Q)_\RR),\qquad
z\longmapsto Ad(h(z))=[M\mapsto h(z)Mh(z)^{-1}].
$$
Since $h(t)=t^2I_V$ for $t\in\RR^*$,
the map $Ad(h(t))$ is the identity, hence  $so(Q))$ is a Hodge structure
of weight $0$. 

The spaces $V^{2,0}$ and $V^{0,2}$ are isotropic for $Q$ and $Q$ gives a duality between them, see \ref{lempol}. So we can choose a basis 
$f_1,\ldots, f_m$ of  $V^{2,0}$ and $f_{m+1},\ldots, f_{2m}$ of
$V^{0,2}$  such that $Q(y,y)=
y_1y_{m+1}+y_2y_{m+2}+\ldots+y_my_{2m}$ with $y=\sum y_if_i$.
In case $m>1$ the following linear map $B$ lies in $so(Q)_\CC$:
$$
B:V_\CC\longrightarrow V_\CC,\qquad B(V^{1,1})=0,\quad B(f_1)=f_{2m},\quad
B(f_m)=-f_{m+1},\quad B(f_i)=0
$$
if $i\neq 1,\,m$. Since $h(z)$ acts as $z^p\bar{z}^q$ on $V^{p,q}$
it is easy to verify that 
$$
h(z)Bh(z)^{-1}=z^{-2}\bar{z}^2B,\qquad{\rm hence}\quad
so(Q)^{-2,2}\neq 0.
$$

Let $W$ be a Hodge structure of weight 1
and let
$h_W:\CC^*\rightarrow GL(W_\RR)$ be the homomorphism 
which defines the Hodge structure on $W$.
The Hodge 
structure on $W\otimes W$ is defined by 
$h_2(z)(u\otimes w)=(h_W(z)u)\otimes(h_W(z)w)$.
Thus $MT(W\otimes W)\subset GL(W)$ where $GL(W)$ acts on $W\otimes W$ by
$A\cdot (u\otimes w)=(Au)\otimes (Aw)$. Therefore the Lie algebra $Lie(GL(W))=End(W)$ acts 
on $W\otimes W$ by
$$
X\cdot(u\otimes w)=(Xu)\otimes w\,+\,u\otimes (Xw)\qquad(X\in End(W),\;
u,\;w\in W).
$$
  
For $X\in End(W)$ the map $Ad(h_2(z))(X)\;(\in End(W\otimes W))$ is given by:
$$
(h_2(z)\cdot X\cdot h_2(z)^{-1})\cdot (u\otimes w)=(h_W(z)Xh_W(z)^{-1}u)\otimes w+
u\otimes(h_W(z)Xh_W(z)^{-1}w),
$$
hence $Ad(h_2(z))(X)=
Ad(h_W(z))(X)$ ($=h_W(z)Xh_W(z)^{-1}\in End(W))$
acting on $W\otimes W$. 
The eigenvalues of $h_W(z)$ on $W_\CC=W^{1,0}\oplus W^{0,1}$
are $z$ and $\bar{z}$, so
the eigenvalues
of the map $End(W)_\CC\rightarrow End(W)_\CC$,
$X\mapsto h_W(z)Xh_W(z)^{-1}$
are $z\bar{z}^{-1},\;1,\;z^{-1}\bar{z}$.
Thus the Hodge structure on $End(W)\;(\subset End(W\otimes W))$
has trivial $(-2,2)$-part.
Therefore $so(Q)$ cannot be contained in $End(W)$.
\qed

\subsection{}
We consider the polarized weight two Hodge structures $(V,h,Q)$
with $\dim V^{2,0}=1$ in some more detail. The quadratic form defined
by the polarization on $V_\RR$ is negative definite on $V_2$, a two dimensional subspace, and positive definite on $V_1$ (cf.\ \ref{lempol}). 

\subsection{Lemma.}\label{equiv}
Let $(V,Q)$ be a $\QQ$-vector space of dimension $n$ with a bilinear
form $Q$ of signature $(2-,(n-2)+)$.
Then there is a natural bijection between the following two sets
\begin{enumerate}
\item
The set of algebraic homomorphisms $h:\CC^*\rightarrow GO(Q)(\RR)$
such that $(V,h,Q)$ is a polarized Hodge structure of weight two.
\item
The set of oriented two-dimensional subspaces $W\subset V_\RR$ such that
the restriction of $Q$ to $W$ is negative definite.
\end{enumerate}

\ts
Given $h$, one defines $W:=V_2$ and the orientation on $V_2$ is defined by 
the basis $v,\,h(e^{\pi i/4})v$, for (any) $v\in V_2, \;v\neq 0$.

Conversely, given the oriented $W$, define $h(e^{i\phi})$ to be rotation with angle $2\phi$ (in the positive sense) on $W$ and to be the identity on $W^\perp$. For $t\in \RR$, 
let $h(t)$ be scalar multiplication by $t^2$ on $V$, this defines
the representation $h:\CC^*\rightarrow GL(V_\RR)$.
\qed

\subsection{Existence of polarized weight two Hodge structures.}\label{exis2}
The Lemma implies that given a $\QQ$-vector space $V$ of dimension $n$ and a
non-degenerate quadratic form $Q$ on $V$ with signature $(2-,(n-2)+)$
there exist polarized weight two Hodge structures $(V,h,Q)$ since we can
certainly find a $W\subset V_\RR$ on which $Q$ is negative definite.

There is a basis of $V_\RR$
such that $Q=-X_1^2-X_2^2+X_3^2+\ldots+X_n^2$. If the restriction of $Q$ to
a subspace $W$ is negative definite then we can find a
basis $a,\,b$ of $W$ such that $a=(1,0,a_3,\ldots a_n)=(1,0,a')$,
$b=(0,1,b_3,\ldots,b_n)=(0,1,b')$
(else $W$ contains a non-zero element $(0,0,c_3,\ldots,c_n)$ contradicting
that $Q$ is negative definite on $W$). 
Thus $W$ depends on $2(n-2)$- real parameters
$(a',b')\in U\subset \RR^{2(n-2)}$ for a certain open subset $U$.

\subsection{Remark.}\label{sotran}
The group $GO(Q)(\RR)$ acts in a natural way on both of the sets mentioned
in Lemma \ref{equiv}:
$$
h\longmapsto h^g:=[z\mapsto gh(z)g^{-1}],\qquad W\longmapsto gW,
$$
this action is compatible with the bijection we indicated. By Witt's theorem,
for any two subspaces $W,\,W'\in V_\RR$ on which $Q$ is negative definite, 
there exists a $g\in SO(Q)(\RR)$ with $gW=W'$ (\cite{Lam}, I.4.2). 
Moreover, if $\dim V>2$ we can also find a $g\in SO(Q)(\RR)$ 
with $gW=W$ and which reverses the orientation on $W$.
Hence $SO(Q)$ acts transitively on the set of Hodge structures
and this set is thus identified with $SO(Q)(\RR)/S(O(2)\times O(n-2))(\RR)$
where $n=\dim V >2$. This set is actually the disjoint union (corresponding to
the choice of orientation) of two copies of a Hermitian symmetric domain,
in particular, the $W$ depend on $n-2$ complex parameters.

\subsection{Lemma.} Let $(V,h,Q)$ be a weight two Hodge structure with
$\dim V^{2,0}=1$. Then, for general $h$, we have $MT(V)=GO(Q)$.

\ts 
The proof is similar to the proof of \cite{vG}, 6.11.
\qed

\section{From weight two to weight one.}

\subsection{}
We recall the construction of Kuga and Satake
(\cite{KS}, \cite{Sata} and \cite{DK3})  which associates to
a polarized Hodge structure $(V,h,Q)$ of weight two,
with $\dim V^{2,0}=1$, a polarized Hodge structure
$(C^+(Q),h_s,E)$ of weight one. 
In Proposition \ref{mtspin} we show that one can recover the
Hodge structure on 
$V$ from the one on $C^+(Q)$ and  that
$V\subset C^+(Q)\otimes C^+(Q)$ (inclusion of Hodge structures).

\subsection{Quadratic Forms.}
Let $(V,h,Q)$ be a polarized weight two rational Hodge structure
with $\dim V^{2,0}=1$, let $n=\dim V$. 
We will simply write $Q(v)$ for $Q(v,v)$, thus the polarization
$Q$ is also viewed as a
(non-degenerate) quadratic form $Q$ on the 
$\QQ$-vector space $V$.
Since $Q$ has signature $(2-,(n-2)+)$ there is a basis of $V$ such that
(\cite{Lam}, Cor.\ 2.4): 
$$
Q:\quad  d_1X_1^2+d_2X_2^2+\ldots d_nX_n^2,\qquad 
d_1,\;d_2<0,\quad d_3,\ldots,d_n>0.
$$

\subsection{Clifford algebras.}
Associated to $(V,Q)$ there is a $2^n$-dimensional associative
$\QQ$-algebra, the Clifford 
algebra $C(Q)$ (cf.\ \cite{Lam}, Ch.\ V, \cite{Sch}, Ch.\ 9) and a linear
injective map $i:V\hookrightarrow C(Q)$.
This is characterized by the following universal property:

Let $A$ be a $\QQ$-algebra and let $f:V\rightarrow A$ be a linear map such that
$f(v)^2=Q(v)$ for all $v\in V$. Then there is a unique $\QQ$-algebra 
homomorphism
$g:C(Q)\rightarrow A$ such that $f=g\circ i$.

 We observe that if $V$ has 
basis $e_1,\ldots,e_n$ and $Q$ is given by $\sum d_iX_i^2$ with respect to 
this basis, one has $e_i^2=Q(e_i)=d_i$ and, for $i\neq j$, $(e_i+e_j)^2=Q(e_i+e_j)=d_i+d_j$
whence:
$$
e_i^2=d_i,\qquad e_ie_j=-e_je_i\qquad{\rm if}\quad i\neq j.
$$
A $\QQ$-basis for $C(Q)$ is given by the products
$$
e^a:=e_1^{a_1}e_2^{a_2}\ldots e_n^{a_n},\qquad a=(a_1,\ldots,a_n)\in\{0,1\}^n
$$
and
$$
C(Q):=\oplus_a\,\QQ e^a.
$$
In particular, the vector space $V$ is a subspace of $C(V)$:
$$
i:V\hookrightarrow C(V),\qquad e_i\longmapsto 
e_1^{0}\ldots e_i^1\ldots e_n^{0},
$$
note that $V$ is not a subalgebra however.
The even Clifford algebra is the sub-algebra:
$$
C^+(Q):=\oplus_{a}\,
\QQ e^a,\qquad a\in\{0,1\}^n,
\qquad \sum a_i\equiv\,0\;{\rm mod}\;2.
$$
Similarly, the vector space $C^-(Q)$ is the span of the
$e^a$ with $\sum a_i\equiv 1$ mod 2.
We will study the Clifford algebra in more detail in section \ref{clif}.

\subsection{The complex structure on $C^+(Q)_\RR$.}
Let 
$$
h:\CC^*\longrightarrow GO(Q)(\RR)\subset GL(V_\RR)
$$ 
be the homomorphism defining
the Hodge structure on $V$.  Recall that 
$$
V_\RR=V_2\oplus V_1,\qquad{\rm with}\quad
V_1\otimes\CC=V^{1,1},\quad
V_2\otimes \CC=V^{2,0}\oplus V^{0,2}, 
$$ 
and the direct sum is orthogonal for $Q$. 
The space $V_2$ is a real two dimensional vector space.

\subsection{Lemma.}\label{defJ}
Let $\{f_1,\,f_2\}$ be a basis of $V_2$ such that
$V^{2,0}=\langle f_1+if_2\rangle$ and $Q(f_1)=-1$. 
\begin{enumerate}
\item
For all  $x,\, y\in\RR$: $Q(xf_1+yf_2)=-(x^2+y^2)$. 
\item The
element $J:=f_1f_2=-f_2f_1\in C^+(Q)_\RR$ satisfies $J^2=-1$.
\item
The element $J$ does not depend on the choice of $f_1,\,f_2$ as above.
\end{enumerate}

\ts Let $\langle v+iw\rangle$ be any basis of $V^{2,0}$ with $v,\,w\in
V_2$.  The polarization restricted to $V_2$ is negative definite,
hence we can find a $\lambda\in\RR$ with $Q(\lambda
v)=\lambda^2Q(v)=-1$, take $f_1=\lambda v,\;f_2=\lambda w$.
From Lemma \ref{lempol} we have that $Q=0$ on $V^{2,0}$, hence 
$0=Q(f_1+if_2,f_1+if_2)=Q(f_1)-Q(f_2)+2iQ(f_1,f_2)$,
hence also $Q(f_2)=-1$ and $f_1,\,f_2$ are perpendicular. It follows that
$Q(xf_1+yf_2)=-(x^2+y^2)$.
Thus in $C(Q)_\RR$ we have: $(xf_1+yf_2)^2=-(x^2+y^2)$
so $f_1f_2+f_2f_1=0$.

In $C(Q)_\RR$ we now have
$f_i^2=-1$ and $f_1f_2=-f_2f_1$, therefore 
$$
J^2=(f_1f_2)(f_1f_2)=-f_1^2f_2^2=-1.
$$
Let $f'_1,\,f'_2$ be another such basis of $V_2$, then there is an orthogonal 
$2\times 2$ matrix $A$ with $f_i'=Af_i$. Since both $f_1+if_2$ and
$f_1'+if_2'$ span $V^{2,0}$, an eigenspace for all $h(z)$, $A$ must commute
with all $h(z)$'s, hence $A$ is a rotation.
Therefore 
$f_1'=af_1+bf_2$ and $f_2'=-bf_1+af_2$ with $a^2+b^2=1$. Thus $f_1'f_2'=
-ab(-1)+(a^2+b^2)f_1f_2+ab(-1)=f_1f_2$. \qed

\subsection{The weight one Hodge structure on $C^+(Q)$.}
With $J$ as in Lemma \ref{defJ} we define a homomorphism
$$
h_s:\CC^*\longrightarrow GL(C^+(Q)_\RR),\qquad
a+bi\longmapsto a-bJ:=[x\longmapsto (a-bJ)x],
$$
(with $a,\,b\in\RR,\;x\in C^+(Q)_\RR$; the `$-$' sign is needed for \ref{mtspin}).
So we let $a-bJ\in C^+(Q)_\RR$ 
act by right multiplication on $C^+(Q)_\RR$.
This is obviously an algebraic homomorphism and defines a
rational Hodge structure of
weight one on the $\QQ$-vector space $C^+(Q)$.

\subsection{Polarizations.}
We show that the Hodge structure $(C^+(Q),h_s)$ has a polarization.
Given $c\in C^+(Q)$ the right multiplication by $c$ is a
$\QQ$-linear map $C^+(Q)\rightarrow C^+(Q)$, $x\mapsto cx$.
We denote by $Tr(c)\;(\in\QQ)$
the trace of this linear map then 
$$
Tr:C^+(Q)\longrightarrow \QQ
$$ 
is a $\QQ$-linear map. 
There is a $\QQ$-linear algebra anti-involution $\iota$ on $C^+(V)$ (so
$\iota(xy)=\iota(y)\iota(x)$) 
which is given by (\cite{Lam}, V, 1.11; \cite{Sch}, 9.3):
$$
\iota:C^+(Q)\longrightarrow C^+(Q),\qquad 
e_1^{a_1}\ldots e_{n}^{a_n}\longmapsto e_n^{a_n}\ldots e_1^{a_1}
\qquad(a_i\in\{0,1\}).
$$

\subsection{Lemma.} \label{lemtrace}
Let $\{e^a\}$ with $a\in\{0,1\}^n,\;\sum a_i\equiv 0$ mod $2$, be the standard basis
of $C^+(Q)$.
\begin{enumerate}
\item We have:
$$
Tr(e^a)=\left\{\begin{array}{ccl} 0&\;{\rm if}&a\neq 0,\\
                                  2^{n-1}&\;{\rm if}&a=0.\end{array}\right.
$$
In particular: 
$$
Tr(xy)=Tr(yx)\qquad{\rm and}\qquad Tr(\iota(x))=Tr(x).
$$
\item For $a,\,b\in\{0,1\}^n$ we have:
$$
Tr(\iota(e^a)e^b)=\left\{\begin{array}{ccl} 0&\;{\rm if}&a\neq b,\\
                                 2^{n-1} d_1^{a_1}d_2^{a_2}\ldots d_n^{a_n}&\;{\rm if}&a=b.\end{array}\right.
$$
\end{enumerate}

\ts
Consider the matrix of multiplication by $e^a$ w.r.t.\ the standard basis.
Since 
$$
e^ae^b=\lambda e^c,\qquad{\rm with}\quad c_i\equiv a_i+b_i\;{\rm mod}\;2
$$
and $\lambda$ is, up to sign, the product of the $d_i$ for which $a_i=b_i=1$,
we see that $e^ae^b$ is a scalar multiple of the basis vector $e^c$ but
$e^c=e^b$ only if $a=0$. Thus $Tr(e^a)=0$ unless $a=0$ and then
$e^0=1$ hence $Tr(e^0)=\dim C^+(Q)$. Since $Tr(e^ae^b)=0=Tr(e^be^a)$
if $a\neq b$ and, obviously, $Tr(e^ae^a)=Tr(e^ae^a)$, we get $Tr(xy)=Tr(yx)$
as $(x,y)\mapsto Tr(xy)$ is bilinear. Finally $\iota(e^a)=\pm e^a$,
hence $Tr(\iota(e^a))=0=Tr(e^a)$ if $a\neq 0$ and if $a=0$ we have $e^a=1=\iota(e^a)$.

For the second part, since $\iota(e^a)=\pm e^a$ we have $\iota(e^a)e^b=\mu e^c$ for some $\mu\in \QQ$ and $c_i=a_i+b_i$ mod $2$. Thus $Tr(\iota(e^a)e^b)=0
$ unless $a=b$. In that case
$\iota(e^a)e^a=d_1^{a_1}\ldots d_n^{a_n}$ since:
$$
\iota(e_1^{a_1}\ldots e_{n}^{a_n})e_1^{a_1}\ldots e_{n}^{a_n}=
e_n^{a_n}\ldots e_1^{a_1}e_1^{a_1}\ldots e_{n}^{a_n}=
e_n^{a_n}\ldots e_2^{a_2}d_1^{a_1}e_2^{a_2}\ldots e_{n}^{a_n}=\ldots=
d_1^{a_1}\ldots d_n^{a_n}.
$$
This concludes the proof of the lemma.
\qed

\subsection{Proposition.} \label{pol}
Let $\alpha:=\pm e_1e_2\in C^+(\QQ)$. Then 
the bilinear form:
$$
E:C^+(\QQ)\times C^+(\QQ)\longrightarrow \QQ,
\qquad
E(v,w):=Tr(\alpha\iota(v)w)
$$
is a polarization of the weight one Hodge structure $(C^+(\QQ),h_s)$
(for suitable choice of sign).

\ts
The form $E$ is obviously $\QQ$-bilinear.
We have $J=f_1f_2=-f_2f_1\in C^+(Q)_\RR$ hence
$\iota(J)=-J$. Thus, with $z=a+bi\in\CC^*$ we get:
$$
\begin{array}{rcl}
E(h_s(z)x,h_s(z)y)&=&Tr(\alpha\iota((a-bJ)x)(a-bJ)y)\\
&=& Tr(\alpha\iota(x)(a+bJ)(a-bJ)y)\\
&=&Tr(\alpha\iota(x)(a^2+b^2)y)\\
&=&z\bar{z}E(x,y).
\end{array}
$$
The symmetry of the bilinear form $E(x,h_s(i)y)$ follows from:
$$
\begin{array}{rcl}
E(x,h_s(i)y)&=&-Tr(\alpha\iota(x)Jy)\\
&=&-Tr(\iota(\alpha\iota(x)Jy))\\
&=&-Tr(\iota(y)\iota(J)\iota^2(x)\iota(\alpha))\\
&=&-Tr(\iota(y)(-J)x(-\alpha))\\
&=&Tr(\alpha\iota(y)h_s(i)x)\\
&=&E(y,h_s(i)x),
\end{array}
$$
here we used $\iota(\alpha)=\iota(e_1e_2)=e_2e_1=-e_1e_2=-\alpha$ and
Lemma \ref{lemtrace}.

It remains to check that
$E(x,h(i)x)$ is either positive or negative definite.
First we show that, given $(V,Q)$, it suffices to consider just one 
Hodge structure $(V,h,Q)$. Let $(V,h,Q)$ and $(V,h',Q)$ be polarized
Hodge structures and let 
$J,\,J'\in C^+(Q)_\RR$ be the associated complex structures on $C^+(Q)_\RR$.
From \ref{sotran} we have a
$g\in SO(Q)(\RR)$ with $h'=ghg^{-1}$ and $gV_2=V_2'$,
in particular the $f_1,\,f_2\in V_2$ with $J=-f_1f_2\in C^+(Q)_\RR$
are mapped to $f_i'=gf_i$ and $J'=-f_1'f_2'$.
There is a $\tilde{g}\in C^+(Q)_\RR$ with $\tilde{g}v\tilde{g}^{-1}=gv$
for all $v\in V\;(\subset C(Q))$ (\cite{Ca}, Theorem 10.3.1)
and  $\tilde{g}\iota(\tilde{g})\in\RR^*$ (\cite{Ca}, Cor.\ 3.1.1,
\cite{Sch}, Lemma 9.3.2). Therefore $\tilde{g}^{-1}=\lambda \tilde{g}$
for a $\lambda \in \RR^*$ and $J'=\tilde{g}J\tilde{g}^{-1}=\lambda\tilde{g}J\iota(\tilde{g})$
so:
$$
\begin{array}{rcl}
E(x,J'x)&=&Tr(\alpha \iota(x)J'x)\\
&=& \lambda Tr(\alpha \iota(x)\tilde{g}J\iota(\tilde{g})x)\\
&=& \lambda Tr(\alpha \iota(y)Jy)\qquad\qquad
{\rm with}\qquad y:=\iota(\tilde{g})x\\
&=& \lambda E(y,Jy)
\end{array}
$$
hence $E(x,J'x)$ is definite iff $E(x,Jx)$ is definite.

We consider the Hodge structure with $V_\RR=V_2\oplus V_1$,
$V_2=\langle e_1,\,e_2\rangle_\RR$ (so $Q$ restricited to $V_2$ is
$d_1X_1^2+d_2X_2^2$ with $d_1,\,d_2<0$) and we take $J=ce_1e_2$ with
$c=(d_1d_2)^{-1/2}\;(\in\RR_{>0})$. Then we have:
$$
E(x,h(i)y)=-Tr(\alpha\iota(x)Jy)=-cTr(e_1e_2\iota(x)e_1e_2y).
$$

We write the basis elements as:
$$
e^a=e_1^{a_1}e_2^{a_2}f_a\qquad{\rm with}\quad f_a=e_3^{a_3}\ldots e_n^{a_n}.
$$
Then $\iota(e^a)=\iota(f_a)e_2^{a_2}e_1^{a_1}$. In case $a_1=a_2=0$ we have:
$$
\iota(e^a)e_1e_2=\iota(f_a)e_1e_2=(-1)^re_1\iota(f_a)e_2=(-1)^{2r}\iota(f_a)=
e_1e_2\iota(e^a),
$$
here $r=a_3+\ldots+a_n$. Similarly, if $a_1=1,\,a_2=0$:
$$
\iota(e^a)e_1e_2=\iota(f_a)e_1e_1e_2=(-1)^re_1\iota(f_a)e_1e_2=
-(-1)^re_1\iota(f_a)e_2e_1=-e_1e_2\iota(e^a),
$$
and proceeding in this way one verifies that:
$\iota(e^a)e_1e_2=
e_1e_2(-1)^{a_1+a_2}\iota(e^a)$.
Therefore:
$$
\begin{array}{rcl}
E(e^a,h_s(i)e^b)&=&-cTr(e_1e_2\iota(e^a)e_1e_2e^b)\\
&=&
-cTr((e_1e_2)^2(-1)^{a_1+a_2}\iota(e^a)e^b)\\&=&
+cd_1d_2(-1)^{a_1+a_2}Tr(\iota(e^a)e^b).
\end{array}
$$
Using the previous lemma we get:
$$
E(e^a,h(i)e^b)=
\left\{\begin{array}{ccl} 0&\;{\rm if}&a\neq b,\\                               2^{n-1}(cd_1d_2)((-1)^{a_1+a_2}d_1^{a_1}d_2^{a_2})d_3^{a_3}\ldots d_n^{a_n}
&\;{\rm if}&a=b.\end{array}\right.
$$
Thus $E(e^a,h(i)e^a)>0$ for all $a$ since $d_1,\,d_2<0,\;d_3,\ldots,d_n>0$. This, combined with
the $\QQ$-bilinearity of $Tr$, proves that $E(x,h(i)x)>0$ for all 
$x\in C^+(Q)-\{0\}$.
\qed

\section{The Mumford-Tate group of the Kuga-Satake Hodge structure.}

\subsection{} In the previous section we constructed a polarized rational
weight one Hodge structure $(C^+(Q),h_s,E)$. 
We recall some basic facts on this Hodge structure.
A detailed study of the Clifford algebra in the next section 
will give more precise information on the simple 
sub-Hodge structures of $C^+(Q)$.

\subsection{The spin representation.}
Let $Q$ a non-degenerate quadratic form on a $F$-vector space $V$
(we usually consider the case $F=\QQ$ but the cases $F=\RR,\,\CC$ are also
of interest to us).
The Spin group $CSpin(Q)$, an algebraic group defined over $F$,
 can be defined as the subgroup of $C^+(Q)^*$,
the units of the ring $C^+(Q)$:
$$
CSpin(Q)=\{g\in C^+(Q)^*:\; gVg^{-1}\subset V\}.
$$
By its very definition, we have a homomorphism
$$
\rho:CSpin(Q)\longrightarrow GL(V),\qquad g\longmapsto[v\mapsto gvg^{-1}].
$$
Since $Q(\rho(g)v)=(\rho(g)v)^2=gv^2g^{-1}=gQ(v)g^{-1}=Q(v)$,
the image of $\rho$ lies in $O(Q)$, the orthogonal group of $Q$,
one actually has: $\rho(CSpin(Q))=SO(Q)$, (\cite{Ca}, Theorem 10.3.1). 

The group $CSpin(Q)$
also acts by multiplication on the left on $C^+(Q)$, 
this gives a  homomorphism, called
the spin representation (cf.\ \cite{Sch}, 9.3):
$$
\sigma:CSpin(Q)\longrightarrow GL(C^+(Q)),\qquad 
g\longmapsto [x\longmapsto gx].
$$ 
We will identify $CSpin(Q)$ with its image in $GL(C^+(Q))$.

\subsection{Proposition.}\label{mtspin}
\begin{enumerate}
\item
  The image of the homomorphism $h_s:\CC^*\rightarrow GL(C^+(Q)_\RR)$ 
is contained in the algebraic group $CSpin(Q)(\RR)$.
Hence $MT(C^+(Q))\subseteq CSpin(Q)$, and 
$$
MT(C^+(Q))= CSpin(Q) \qquad{\rm if}\quad MT(V)=GO(Q).
$$
\item\label{rho}
For $t,\,\phi\in\RR$ we have:
$$
h(te^{i\phi})=t^2\rho(h_s(e^{i\phi})),
$$
so we recover $(V,Q)$ from the Hodge structure on $C^+(Q)$ defined by
$h_s$.
\item
The Hodge structure $(V,h)$ is a sub-Hodge structure of 
$(C^+(Q)\otimes C^+(Q),h_s\otimes h_s)$:
$$
V\hookrightarrow C^+(Q)\otimes C^+(Q).
$$
\end{enumerate}

\ts First we show $h_s(\CC^*)\in CSpin(Q)(\RR)$.
Since $h_s(te^{i\phi})=t^2(a-bJ)$ with $a^2+b^2=1$ and then
$(a-bJ)^{-1}=a+bJ$ it suffices to show $(a-bJ)V(a+bJ)\subset V$.
For $v\in V_\RR\subset C(Q)_\RR$ and $a,\,b\in\RR$ we have:
$$
(a-bJ)v(a+bJ)=(a^2v-b^2JvJ)+ab(-Jv+vJ).
$$
Recall that $V_\RR=V_1\oplus V_2$ (orthogonal sum) and $J=f_1f_2$ with
$f_i\in V_2$. Thus $Jv=vJ$ for $v\in V_1$, and hence 
$(a-bJ)v(a+bJ)=(a^2+b^2)v\in V$. Note that 
$$
Jf_1=(f_1f_2)f_1=-f_1^2f_2=f_2\in V,\qquad
Jf_2=(f_1f_2)f_2=-f_1\in V.
$$
Similarly, $f_1J=-f_2,\;f_2J=f_1$. This gives:
$$
(a-bJ)f_1(a+bJ)=(a^2-b^2)f_1-2abf_2,\qquad
(a-bJ)f_2(a+bJ)=(a^2-b^2)f_2+2abf_1
$$
hence also $(a-bJ)V_2(a+bJ)\subset V$ 
and we conclude that $h_s(\CC^*)\in CSpin(Q)(\RR)$.

Moreover, we see that  $\rho(h_s(a+bi))v=v$ for $v\in V_1$ and that
 $f_1\pm if_2$
is an eigenvector of $\rho(h_s(a+bi))\;(\in GL(V_\RR))$ with eigenvalue
$(a\pm ib)^2$, this verifies \ref{mtspin}.\ref{rho}.

We already saw that $h_s(\CC^*)\subset CSpin(Q)(\RR)$, hence
$MT(C^+(Q))\subseteq CSpin(Q)$.
If we assume $MT(V)=GO(Q)$ and
$G$ is a subgroup of $SO(Q)$
defined over $\QQ$ with $h(e^{i\phi})\in G(\RR)$
for all $\phi\in \RR$, then $G$  must be equal to $SO(Q)$.
Hence $\rho(MT(C^+(Q)))=SO(Q)$. Since 
$\ker(\rho)\cong\GG_m\;(\subset MT(C^+(Q)))$,
we get $MT(C^+(Q))= CSpin(Q)$.

For $g\in CSpin(Q)$ one has $\nu(g):=\iota(g)g\in\GG_m$ 
(\cite{Sch}, Lemma 9.3.2), hence:
$$
E(gv,gw)=Tr(\alpha\iota(v)\iota(g)gw)=\nu(g)E(v,w).
$$ 
Therefore the isomorphism $V\rightarrow V^*$ defined by $E$ is
equivariant for the action of $CSpin(Q)$ (up to the homomorphism $\nu$).
We will identify the $CSpin(Q)$-representations
$End(C^+(Q))=C^+(Q)^*\otimes C^+(Q)\cong C^+(Q)\otimes C^+(Q)$.

We choose an invertible element, say $e_1$, in 
$V\;(\subset C(Q))$. Then we have an inclusion:
$$
V\hookrightarrow End(C^+(Q)),\qquad 
v\longmapsto M_v:=[y\mapsto vye_1].
$$
The image of $V$ in $End(C^+(Q))$ is a sub-representation
on which $CSpin(Q)$ acts via $\rho$:
$$
(gM_vg^{-1})(y)=(gM_v)(g^{-1}y)=g(vg^{-1}ye_1)=(gvg^{-1})ye_1=M_{gvg^{-1}}(y)
=M_{\rho(g)v}y.
$$
As $h_s(\CC^*)\subset CSpin(Q)(\RR)$ it follows that
$V\hookrightarrow C^+(Q)\otimes C^+(Q)$ is sub-Hodge structure.
\qed

\subsection{}
The spin representation is not irreducible in general. 
In fact, for any $x,\,y\in C^+(Q)$ we have $\sigma(g)(xy)=gxy=(\sigma(g)x)y$,
so the $\QQ$-linear maps $C^+(Q)\rightarrow C^+(Q)$, $x\mapsto xy$ commute
with $CSpin(Q)$. Therefore we have an injective map 
$$
C^+(Q)\hookrightarrow End_{CSpin}(C^+(Q))\quad
(:=\{M\in End(C^+(Q)):\;M\sigma(g)=\sigma(g)M\;\}).
$$
(Due to the action on the right, we should write $C^+(Q)^{op}$,
but $C^+(Q)^{op}\cong C^+(Q)$, cf.\ \cite{Lam}, V, Prop.\ 1.11.)
These are all the maps which commute with the $CSpin(Q)$
representation:

\subsection{Lemma.}
We have:
$$
C^+(Q)\cong End_{CSpin(Q)}(C^+(Q)).
$$

\ts
The group $CSpin(Q)$ is an algebraic subgroup of $GL(C^+(Q))$ and thus
its Lie algebra, $cspin(Q)$ is a subalgebra of $End(C^+(Q))$ and
$cspin(Q)\otimes_\QQ\CC$ is the Lie algebra of the complex Lie group $
CSpin(Q)(\CC)$. Moreover,
$End_{CSpin(Q)}(C^+(Q))\cong End_{cspin(Q)}(C^+(Q))$.
The latter is the subspace of $C^+(Q)$ defined by:
$$
End_{cspin(Q)}(C^+(Q))=
\{X\in End(C^+(Q)):\; XM-MX=0\quad \forall M\in cspin(Q)\}.
$$
Considering these equations for $X\in End(C^+(Q)_\CC)$ we see that
$$
End_{cspin(Q)}(C^+(Q))\otimes_\QQ\CC\cong End_{cspin(Q)\otimes_\QQ\CC}
(C^+(Q)\otimes_\QQ\CC).
$$
From the representation theory of complex Lie algebras
(see \ref{clifsplit} below for the case $n$ even) we know that 
$$
End_{cspin(Q)\otimes_\QQ\CC} (C^+(Q)\otimes_\QQ\CC)
\cong C^+(Q)\otimes_\QQ\CC,
$$
we use that $C^+(Q)\otimes_\QQ\CC$ is the even Clifford algebra of the
quadratic form $Q$ on $V\otimes_\QQ\CC$. Thus $\dim_\QQ 
End_{cspin(Q)}(C^+(Q))=\dim_\CC C^+(Q)\otimes_\QQ\CC=2^{n-1}$, hence
$C^+(Q)$ must be all of $End_{cspin(Q)}(C^+(Q))$ for dimension reasons.
\qed

\subsection{Example.}\label{clifsplit}
In case $Q$ is 
the quadratic form $Y_1Y_{m+1}+\ldots+Y_mY_{2m}$ (this is the case over a finite extension of $\QQ$), the Clifford algebra and the spin representation
can be described as follows
(cf.\ \cite{FH}, p.\ 304;
they consider the 
Clifford algebra over $\CC$, but the same arguments work over any extension
of $\QQ$).

Let  $Z$ be the $m$-dimensional subspace of
$V$ defined by $Y_{m+1}=\ldots=Y_{2m}=0$  Then
$$
C(Q)\cong End(\Lambda^{\ast} Z),\qquad 
\Lambda^{\ast}Z=\QQ\oplus Z\oplus\wedge^2Z\oplus\ldots
$$
The even Clifford algebra is identified with the subalgebra
$$
C^+(Q)\cong End(\Lambda^{even}Z)\times End(\Lambda^{odd}Z)\cong 
M_{2^{m-1}}(\CC)\times M_{2^{m-1}}(\CC).
$$
Since $C^+(Q)\subset End_{cspin(Q)}(C^+(Q))$, 
the spin representation must, at least, split in
the direct sum of $2^{m-1}$ copies of a representation $\sigma_+$ and
$2^{m-1}$ copies of a representation $\sigma_-$. By considering the
Lie algebra action one can see that $\sigma_+$ and $\sigma_-$ are both 
irreducible and are not isomorphic (see \cite{FH}). 
This implies that $C^+(Q)= End_{cspin(Q)}(C^+(Q))$.
The representations 
$\sigma_+$ and $\sigma_-$, each of dimension $2^{m-1}$, are the half spin
representations of $so(Q)$.

\section{Clifford algebras}\label{clif}
\subsection{}
We recall that the Clifford algebra can be built up from quaternion algebras.
If $Q= d_1X_1^2+d_2X_2^2+\ldots d_nX_n^2$, with $d_i\in\QQ$,
we write:
$$
V=\langle d_1\rangle\oplus\ldots\oplus \langle d_n\rangle.
$$

\subsection{Quaternion Algebras.}\label{quats}
The quaternion algebra $A=(a,b)_F$
 over a field $F$ has an $F$-basis
$1,\,i,\,j,\,k$ such that
$$
A=F\oplus Fi\oplus Fj\oplus Fk,\qquad i^2=a,\quad j^2=b,\quad ij=k=-ji,
\qquad a,\,b\in F^*,
$$
with $F^*=F-\{0\}$ and $xa=ax$ for all $x\in F,\,a\in A$.
In case $F=\QQ$ we omit the index $F$. 
Note however that $(a,b)_F$ and $(c,d)_F$ may be isomorphic without
$a=c$ and $b=d$. For example,
$(a,b)_F\cong (b,a)_F$ (via $i\mapsto j,\;j\mapsto i,\;k\mapsto -k$).
The algebra $M_2(F)$ of $2\times 2$ matrices over $F$ is a
quaternion algebra and for any $b\in F^*$     
$$
M_2(F)\cong (1,b)_F,\qquad i=\left( \begin{array}{rr} 1 & 0 \\
                                    0 & -1
                  \end{array} \right),\qquad
j=\left( \begin{array}{rr} 0 & 1\\
                                    b & 0 
                  \end{array} \right).
$$
The quaternion algebra $(a,b)_F$ is either a (skew) field or is isomorphic to
to the matrix algebra $M_2(F)$. One has: $(a,b)_F\cong M_2(F)$ if and only if
the equation $ax^2+by^2=abz^2$ has a non-trivial solution $(x,y,z)\in F^3$
(\cite{Lam}, III,2.7; \cite{Sch}, Corollary 2.11.10;
note that if $ax^2+by^2-abz^2=0$ then $(xi+yj+zk)^2=0$ hence $(a,b)_F$
is not a field).

\subsection{Example.} \label{1var}
In case $Q=dX^2$, $d\in\QQ^*$, the Clifford algebra is isomorphic
to $\QQ[T]/(T^2-d)$, which is a field if $d$ is not a square and is isomorphic
to $\QQ\times\QQ$ if $d$ is a square. In this case $C^+(Q)\cong\QQ=M_1(\QQ)$.

\subsection{Example.}\label{2var}
In case $Q=aX_1^2+bX_2^2$,
$C(Q)$ is the quaternion algebra $(a,b)$.
The even Clifford algebra $C^+(Q)=\QQ\oplus\QQ e_1e_2$ with
$(e_1e_2)(e_1e_2)=-e_1^2e_2^2=-ab$, so if $-ab$ is not a square,
the even Clifford algebra is a field: $C^+(Q)\cong\QQ(\sqrt{-ab})$.
In case $-ab=c^2$ is a square, $C^+(Q)\cong\QQ\times\QQ$, $e_1e_2\mapsto(c,-c)$.

\subsection{Example.}\label{S929} Let $Q=d_1X_1^2+d_2X_2^2+ d_3X_3^2$,
note that $1,\;i:=e_1e_2,\;j:=e_2e_3,\;k:=d_2e_1e_3$, satisfy the rules of the quaternion algebra $(-d_1d_2,-d_2d_3)$ and that
$1,\ldots,k$ are a basis
of $C^+(Q)$ hence
$$
C^+(Q)\cong (-d_1d_2,-d_2d_3).
$$
The element $z=e_1e_2e_3$ is in the center of the Clifford algebra $C(Q)$
(since it commutes with the generators $e_1,\,e_2,\,e_3$) and $z^2=-e_1^2e_2^2e_3^2=-d_1d_2d_3$. Moreover, $C(Q)=C^+(Q)\oplus zC^+(Q)$,
hence $C(Q)\cong C^+(Q)\otimes \QQ(\sqrt{-d_1d_2d_3})$. Since $\QQ(\sqrt{-d_1d_2d_3})\cong C(-d_1d_2d_3X^2)$, we have
(\cite{Sch}, Lemma 9.2.9, with some changes in notation):
$$
C( d_1X_1^2+d_2X_2^2+ d_3X_3^2)\cong(-d_1d_2,-d_2d_3)\otimes 
C(\langle -d_1d_2d_3\rangle).
$$

\subsection{Graded tensor products}
To determine the Clifford algebra of the quadratic form $Q=d_1X_1^2+\ldots+d_nX_n^2$ on $V$, so $V=\langle d_1\rangle
\oplus\ldots\oplus \langle d_n\rangle$, we use that if $V=V'\oplus V''$ is the
orthogonal direct sum with $Q_{|V'}=Q'$ and $Q_{|V''}=Q''$ then
(\cite{Sch}, 9.2.5):
$$
C(Q'\oplus Q'')=C(Q')
\hat{\otimes} C(Q'').
$$
Here $\hat{\otimes}$, the graded tensor product,
is the tensor product of the underlying $\QQ$-vector spaces but the product
is given by (\cite{Sch}, 9.1.4))
$(x\hat{\otimes} y)(x'\hat{\otimes} y')=\epsilon xx'\hat{\otimes}yy'$
with $x,\,x' \in C^{\pm}(Q')$ and $y,\,y'\in C^{\pm}(Q'')$ and
$\epsilon=1$ except if $y\in C^-(Q'')$ and $x'\in C^-(Q')$, in that case
$\epsilon=-1$.

An easy example is the case $Q=aX^2+bY^2$, $Q'=\langle a\rangle$ and
$Q''=\langle b\rangle$, in that case $C(Q'\oplus Q'')=
C(Q')\hat{\otimes} C(Q'')\cong (a,b)$ as 
in example \ref{2var}.
Using example \ref{S929} as induction step we then obtain, 
from  Theorem 9.2.10 of \cite{Sch} and a basic fact on quaternion algebras,
the following result:

\subsection{Theorem.} \label{simcomp}
Let $Q=\sum d_iX_i^2$, with $d_i\in\QQ^*$.
Then we have:
\begin{enumerate}
\item 
In case $n=2m$, let $d:=(-1)^md_1\ldots d_n$.
Then the even
Clifford algebra $C^+(Q)$ is isomorphic to one of the following two algebras
\begin{enumerate}
\item
In case $\sqrt{d}\in \QQ$ we have $C^+(Q)\cong M_{2^{m-2}}(D)\times M_{2^{m-2}}(D)$ with a quaternion algebra $D$ over $\QQ$.
\item
In case $\sqrt{d}\not\in \QQ$ we have
$C^+(Q)\cong M_{2^{m-2}}(D)=M_{2^{m-2}}(\QQ)\otimes_\QQ D$, where $D$
is a quaternion
algebra over a quadratic field extension $F=\QQ(\sqrt{d})$ of $\QQ$.
\end{enumerate}
\item 
In case $n=2m+1$, we have
$C^+(Q)\cong M_{2^{m-1}}(D)$
for a quaternion algebra $D$ over $\QQ$.
\end{enumerate}
Note that if $D\cong M_2(\QQ)$ then $M_{2^a}(D)\cong M_{2^{a+1}}(\QQ)$.

\ts
For $n\leq 3$ see \ref{1var}, \ref{2var}, \ref{S929}. For $n>3$ use
Example \ref{S929}:
$$
\begin{array}{rcl}
C(\langle d_1\rangle\oplus\ldots\oplus \langle d_n\rangle)&\cong&
 C(\langle d_1\rangle\oplus\ldots\oplus \langle d_3\rangle)
\hat{\otimes}
C(\langle d_4\rangle\oplus\ldots\oplus \langle d_n\rangle)\\
&\cong&
(-d_1d_2,-d_2d_3)\otimes 
C(\langle -d_1d_2d_3\rangle)\hat{\otimes}
C(\langle d_4\rangle\oplus\ldots\oplus \langle d_n\rangle)\\
&\cong&
(-d_1d_2,-d_2d_3)\otimes 
C(\langle -d_1d_2d_3\rangle\oplus
\langle d_4\rangle\oplus\ldots\oplus \langle d_n\rangle).
\end{array}
$$
in the last lines the first tensor product is the usual one
since all elements in
$(-d_1d_2,-d_2d_3)$ are even elements of 
$C(\langle d_1\rangle\oplus \langle d_2\rangle\oplus \langle d_3\rangle)$.

In case $n=4$ one thus finds $C(Q)\cong (-d_1d_2,-d_2d_3)\otimes
(-d_1d_2d_3,d_4)$ and $C^+(Q)\cong (-d_1d_2,-d_2d_3)\otimes F$, with
$F\cong \QQ[X]/(X^2-d_1d_2d_3d_4)$ the even Clifford algebra of 
$\langle -d_1d_2d_3\rangle\oplus \langle d_4\rangle$. So $F$ is either a field
or is isomorphic to $\QQ\times\QQ$.

We continue this game and we find quaternion
algebras $A_{13},\,A_{15},\ldots A_{1,2k+1}$, such that
$$
C(Q)\cong A_{13}\otimes A_{15}\otimes\ldots A_{1,2k+1}
\otimes
C(\langle (-1)^{k}d_1d_2\ldots d_{2k+1}\rangle
\oplus
\langle d_{2k+2}\rangle\oplus\ldots
\langle d_n\rangle),
$$
and the even Clifford algebra is:
$$
C^+(Q)\cong  A_{13}\otimes A_{15}\otimes\ldots A_{1,2k+1}\otimes
C^+(\langle (-1)^{k}d_1d_2\ldots d_{2k+1}\rangle
\oplus
\langle d_{2k+2}\rangle\oplus\ldots
\langle d_n\rangle).
$$

We conclude that if $n=2m$, 
then $C(Q)$ is the tensor product of $m-1$ quaternion algebras which is
then graded tensored with an $m$-th quaternion algebra. Thus
$C^+(Q)$ is a tensor product of $m-1$ quaternion algebras with an algebra
$\QQ[X]/(X^2-d)$ and it is not hard to see that $d=(-1)^md_1\ldots d_{2m}$.

In case $n=2m+1$, $C(Q)$ is the tensor product of $2m$ quaternion algebras 
(which is $C^+(Q)$) graded tensored by $C(\langle(-1)^md_1\ldots d_{2m+1}
\rangle)$. 

Now we recall that the tensor product $A\otimes_\QQ B$
of two quaternion algebras $A$, $B$  over $\QQ$
is isomorphic with $M_2(D)=M_2(\QQ)\otimes_\QQ D$ 
for some quaternion algebra $D$. 
Let $A=(a,c)$, $B=(b,d)$ then, in $A$, we have
$$
(x_1i+x_2j+x_3k)^2=ax_1^2+cx_2^2+acx_3^2.
$$
The quadratic form in
6 variables $ax_1^2+cx_2^2+acx_3^2-(by_1^2+dy_2^2+bdy_3^2)$ is indefinite
for any choice of signs for $a,\ldots ,d$. 
Hence Meyer's theorem (\cite{serre}, Corollary 4.3.2)
implies that it has a non-trivial zero. 
Thus there are $x\in A$, $y\in B$ with $x^2=y^2=e$, so if 
$A$ and $B$ are fields they have the quadratic field
$K=\QQ(\sqrt{e})$ in common. Then $A$ and $B$ are cyclic algebra's
over $K$: $A=(K/\QQ,r)$, $B=(K/\QQ,s)$ (so $A=K\oplus K\alpha$ with
$\alpha^2=r$
and $\alpha x=\bar{x}\alpha$ where `$\bar{\phantom{z}}$' is the
conjugation on $K$).
In \cite{Lam}, III, 2.11 and \cite{Sch}, Theorem 8.12.7
one finds an explicit isomorphism $A\otimes_\QQ B\cong M_2(D)$
with $D=(K/\QQ,rs)$.
\qed

\section{Kuga-Satake varieties}
\label{KSV}
\subsection{}
The rational, polarized, weight one Hodge structure $(C^+(Q),h_s,E)$
defines an isogeny class of abelian varieties.
Each variety in this isogeny class is called a Kuga-Satake variety for
$(V,h,Q)$. More precisely, consider the dual $W$ of the
$2^{n-2}$-dimensional complex
vector space $C^+(Q)^{1,0}$:
$$
W:=\left(C^+(Q)^{1,0}\right)^*,\qquad {\rm let}\quad
\Gamma\subset C^+(Q)^*
$$
be a free $\ZZ$-module with $\Gamma\otimes_\ZZ\QQ=C^+(Q)^*$, the dual of the $\QQ$-vector space $C^+(Q)$. The
image of $\Gamma$ under the projection
from $(C^+(Q)^{0,1})^*$: $\Gamma\subset C^+(Q)^*_\CC\rightarrow W$
is a lattice in the complex vector space $W$ and the quotient is an abelian
variety $A_\Gamma$ (with Riemann form $E_\Gamma$ defined by $E$), which is
a Kuga-Satake variety of $(V,h,Q)$.
There is a natural isomorphism of  rational, polarized, weight one
Hodge structures:
$$
(H^1(A_\Gamma,\QQ),E_\Gamma)\,\cong \,(C^+(Q),h_s,E).
$$
For an abelian variety $A$ one has
$End(A)\otimes\QQ\cong End_{Hod}(H^1(A,\QQ))$ and if
$$
End(A)\otimes_\ZZ\QQ\cong M_{n_1}(D_1)\times\ldots\times M_{n_d}(D_d)
\qquad{\rm then}\quad
A\sim A_1^{n_1}\times\ldots\times A_d^{n_d}
$$
where the $D_i$ are (skew) fields, $\sim$ means isogeneous and the $A_i$
are simple abelian varieites. Thus Theorem \ref{simcomp}
gives the decomposition
in simple factors of a Kuga-Satake variety in case $MT(V)=GO(Q)$.

\subsection{Example.} Let $(V,h,Q)$ be
a three dimensional rational Hodge structure
of weight two with $\dim V^{2,0}=1$ and $MT(V)=GO(Q)$. Assume that with
 respect to some basis of $V$ we can write $Q=-X_1^2-X_2^2+dX_3^2$.
Then from Example \ref{S929} we know:
$$
C^+(Q)\cong (-1,d).
$$
A Kuga-Satake variety of $(V,h,Q)$ is an abelian surface $X$ with 
$End(X)\cong C^+(Q)$. Hence it is a product of two isogeneous elliptic 
curves if $C^+(Q)$ is not a field, otherwise it is a simple abelian surface.
Note that $C^+(Q)\cong M_2(\QQ)$ iff $-x^2+dy^2=-dz^2$ has a non-trivial
solution in $\QQ^3$. This happens for example if $d=1$. In case $d\equiv 3\;$
mod $4$ we get a field, this is best seen by multiplying the equation 
$-x^2+dy^2=-dz^2$ by $d$ and substituting $y:=d^{-1}y,\;z:=d^{-1}z$ to get
$dx^2=y^2+z^2$. If this had a non-trivial solution, we had one in $\ZZ^3$
(multiply by product of denominators) and after dividing by a power of $2$ at
least one of $x,\,y,\,z$ would be odd. Now consider the equation mod $4$ and
use that a square is either $0$ or $1$ mod 4.

\subsection{Geometric realizations.}\label{real}
In \ref{exis2} we observed that given $(V,Q)$, one can find weight two
Hodge structures $(V,h,Q)$. Since $V\subset C^+(Q)\otimes C^+(Q)$ we find
an inclusion of Hodge structures:
$$
V\hookrightarrow H^1(A_\Gamma,\QQ)\otimes H^1(A_\Gamma,\QQ)\subset
H^2(A_\Gamma\times A_\Gamma,\QQ).
$$
Thus any polarized weight two Hodge structure $(V,h,Q)$ with $\dim V^{2,0}=1$
is a sub-Hodge structure of the cohomology of some algebraic variety.

This is {\em not} true if $\dim V^{2,0}>1$, Griffiths work on
variations of Hodge structures (`Griffiths transversality')
implies that the general polarized Hodge structure of weight two
(and similar results hold for higher weight) with
$\dim H^{2,0}>1$ is not a sub-Hodge structure of
the cohomology of an algebraic variety, see \cite{CKT} and \cite{Ma}.

\section{Abelian varieties of Weil type and Kuga-Satake varieties}
\label{Lom}

\subsection{} The Hodge conjecture for abelian varieties of dimension
at least 4 is still open (see \cite{BG} for a recent overview). For 4
dimensional abelian varieties (abelian 4-folds)
Moonen and Zarhin \cite{MZ}, \cite{MZ2}
proved that if the Hodge conjecture is true for abelian 4-folds of Weil type, then the Hodge conjecture is true for all abelian 4-folds. 

An abelian $2n$-fold $A$ is of Weil type if $End(A)\otimes_\ZZ\QQ$ contains an
imaginary quadratic field
$K=\QQ(\sqrt{-d}),\,d>0$ and the action of the
endomorphism $\sqrt{-d}$ on the tangent space at the
origin of $A$ has eigenvalues $\sqrt{-d},\,-\sqrt{-d}$,
each with multiplicity $n$.
Associated to the pair $(A,K)$ is a `discriminant' $\delta\in
\QQ^*/N(K^*)$ where $N(K^*)$ is the group of the norms of $K^*$,
that is, the elements $a^2+b^2d$, $a,\,b\in\QQ$ (see \cite{vG} for an
introduction to these varieties).  

The Hodge conjecture for the general abelian $2n$-fold of Weil type with $d=3$,
$\delta=1$ is proved by C.\ Schoen in \cite{Sab} for $n=2$. The case $n=3$
is done in \cite{Sdec}. There
he also proves that the Hodge conjecture for the general abelian 4-fold of Weil type with $d=3$ (and any $\delta$) follows from this by specializing the
6-fold to a product of a 4-fold and an abelian surface.

The following result, due to G.\ Lombardo, may also be of some interest. 
The case that $K=\QQ(\sqrt{-1})$ was considered by Paranjape \cite{Pa}
and, from a different point of view, by Matsumoto and others (\cite{M} and
references given there).

\subsection{Theorem.} Let $(A,K=\QQ(\sqrt{-d}))$
be an abelian 4-fold of Weil type 
with discriminant $\delta=1$.
Then $A^4$ is a Kuga-Satake variety of a weight two Hodge structure $(V,h,Q)$
with 
$$
Q=-X_1^2-X_2^2+X_3^2+X_4^2+X_5^2+dX_6^2.
$$
Conversely, if $(V,h,Q)$ is a weight two polarized Hodge structure with $Q$
as above, then the Kuga-Satake variety of $(V,h,Q)$ is isogeneous to
$A^4$ with $A$ an abelian 4-fold of Weil type.

\ts See \cite{Lom}, we only verify that the Kuga-Satake variety of a
general $(V,h,Q)$
is isogeneous to $A^4$ with $A$ an abelian 4-fold with $K=End(A)$.
Using the rules from the proof of
Theorem \ref{simcomp}, one finds:
$$
C(Q)\cong (-1,1)\otimes C(-X^2+X_4^2+X_5^2+dX_6^2)\cong 
M_2(\QQ)\otimes  C(-X^2+X_4^2+X_5^2+dX_6^2),
$$
since $(-1,1)\cong (1,-1)\cong M_2(\QQ)$, cf.\ \ref{quats}.
Next:
$$
C(-X^2+X_4^2+X_5^2+dX_6^2)\cong (1,1)\otimes
\left(C(\langle 1\rangle)\hat{\otimes} C(\langle d\rangle)\right)
\cong M_2(\QQ)\otimes C(X^2+dY^2),
$$
since
$C^+(X^2+dY^2)\cong\QQ(\sqrt{-d})\cong K$ 
we get:
$$
C^+(Q)\cong M_2(\QQ)\otimes M_2(\QQ)\otimes C^+(X^2+dY^2)
\cong M_4(\QQ)\otimes K.
$$ 
Hence for general $(V,h,Q)$ any 
Kuga-Satake variety is isogeneous to 4 copies of a simple abelian variety
$A$, so $\dim A=4$, with $End(A)\cong K$.  \qed

\section{The Kuga-Satake-Hodge conjecture}\label{KSH}

\subsection{}\label{introksh}
Given a polarized Hodge structure $(V,h,Q)$ of weight two with 
$\dim V^{2,0}=1$  there exists an abelian variety $A$,
the Kuga-Satake variety of $(V,h,Q)$, with the property that $V$ is a sub-Hodge
structure of $H^2(A\times A,\QQ)$, cf.\ section \ref{real}. In case there
is another algebraic variety $X$ with $V\hookrightarrow H^2(X,\QQ)$, the
Hodge conjecture predicts (cf.\ section \ref{morcor})
the existence of an algebraic cycle $Z\subset A^2\times X$,
the Kuga-Satake-Deligne correspondence,
which realizes the morphism of Hodge structures 
$$
f:H^2(A\times A,\QQ)
\stackrel{\pi}{\longrightarrow}V
\stackrel{i}{\longrightarrow}
H^2(X,\QQ).
$$
(More generally, one can consider $V\hookrightarrow H^{2+2n}(X,\QQ)(n)$.)

A particular case where this happens is when the algebraic variety
$X$ has $\dim H^{2,0}(X,\QQ)=1$. 
The  space
of the Hodge cycles $B(H^2(X,\QQ))$ a sub-Hodge structure of $H^2(X,\QQ)$,
it is the image of $CH^1(X)_\QQ$ and is
called the Neron-Severi group of $X$ (tensored by $\QQ$).
Let $V$ be the orthogonal
complement in $H^2(X,\QQ)$:
$$
H^2(X,\QQ)=V\oplus NS(X)_\QQ.
$$
Then $V$, with the polarization $Q$ 
induced by the one on $H^2(X,\QQ)$, is of the type we consider.
It has the additional convenient
property that it is simple in the sense that if
$W\subset V$ is a sub-Hodge structure, then $W=0$ or $W=V$. In fact,
$V=W\oplus W^\perp$ and if $W^{2,0}=0$ then 
$W\subset NS(X)_\QQ\cap V=0$, else $(W^\perp)^{2,0}=0$
so $W^\perp=0$ and $W=V$.

Since for an ample divisor $Y$ on a variety $X$ the restriction map
$H^2(X,\QQ)\hookrightarrow H^2(Y,\QQ)$ is injective if $\dim Y\geq 2$
(Lefschetz Theorem), we may assume without loss of generality that $X$ is a
surface. Then the Hodge conjecture leads to:

\subsection{Kuga-Satake-Hodge conjecture.}\label{ksh}
Let $X$ be a smooth surface 
with $\dim H^{2,0}(X)=1$ and let $H^2(X,\QQ)=V\oplus NS(X)_\QQ$.
Let $Q$ be the induced polarization on $V$ and
let $A$ be a Kuga-Satake variety of $(V,h,Q)$.

Then there exist a surface $Z$ and a diagram
$$
\begin{array}{rcl} 
Z&\stackrel{\phi}{\longrightarrow}&A\times A\\
\pi\downarrow&&\\
X&&
\end{array}
\qquad\mbox{such that}\quad \pi_*\phi^*:
H^2(A\times A,\QQ)\longrightarrow H^2(X,\QQ)
$$
induces an isomorphism $V\stackrel{\cong}{\longrightarrow} V$.

\subsection{} We show that conjecture \ref{ksh} is indeed
a special case of the
Hodge conjecture. Let $f:H^2(A^2,\QQ)\rightarrow H^2(X,\QQ)$ be the
morphism of Hodge structures as in \ref{introksh}. Then
$f\in B(H^{2d}(A^2\times X,\QQ))$ with $d=\dim A^2$ (cf.\ \ref{morcor}). 
The Hodge conjecture
implies that $f$ should be the
class of a codimension $d$ cycle $\sum_i a_iZ_i$ on the $d+2$-dimensional variety
$A^2\times X$, hence each $Z_i$ is a surface (an irreducible 2 dimensional
variety). Each $Z_i$ defines a morphism of Hodge structures
$[Z_i]:H^2(A^2,\QQ)\rightarrow H^2(X,\QQ)$, since $V$ is simple the restriction
of each $[Z_i]$ to $V$
is either $0$ or is an isomorphism on its image. As $f$ is an isomorphism, there is at least one $Z_i$ such that $[Z_i]$ induces an isomorphism, take $Z$ to be that $Z_i$.

\subsection{} Assume that $Z$ exists.
It is easy to see that we may replace $Z$ by its desingularization.
The map $\pi:Z\rightarrow X$ must be surjective (because $H^{2,0}(X)\subset \pi_*H^2(Z,\CC)$). Since $H^2(A^2,\QQ)=\Lambda^2H^1(A^2,\QQ)$ and
pull-back is compatible with cup product, we get: 
$$
V\hookrightarrow \pi_* {\rm Image}\left(\Lambda^2H^1(Z,\QQ)\longrightarrow
H^2(Z,\QQ)\right).
$$
Conversely, given a surface $Z$ with a surjective map $\pi:Z\rightarrow X$
having the above property, the albanese map $Z\rightarrow Alb(Z)$
is essentially a map $\phi:Z\rightarrow A^2$ with $\pi_*\phi^*$ as in
\ref{ksh}.

In the example of Paranjape \cite{Pa} and in Example \ref{nori} below
  the surface $Z$ is
  a product of two curves, also in the example of Voisin \cite{Vks}
  it seems to be possible to choose for $Z$ a product of curves.
In all these examples the  surfaces $X$ 
 are K3 surfaces (so $\dim H^2(X,\QQ)=22$, $\dim H^{2,0}(X)=1$),
but their Neron-Severi groups are rather large.

\section{Hodge structures and imaginary quadratic fields.}

\subsection{} The example of C.\ Voisin
\cite{Vks} deals with a polarized weight two
Hodge structure $(V,h,Q)$ with $\dim V^{2,0}=1$ which has
an automorphism $\phi$ of order 3 preserving the polarization:
$Q(v,w)=Q(\phi v,\phi w)$ (it is induced from an
automorphism of a $K3$ surface). The action of $\phi$ on $V$ gives
$V$ the structure of a vector space over $K=\QQ(\sqrt{-3})$.

Voisin asked for the simple factors
of the Kuga-Satake variety
of $(V,h,Q)$. She already proved
that the (isogeny class of an) elliptic
curve $A_0$
with complex multiplication by $K$ and an abelian variety $A_1$ with
$2\dim A_1=\dim V$, $K\subset End(A_1)\otimes\QQ$, are simple components.
The following theorem (an exercise in representation theory)
gives the simple factors in general.

\subsection{Theorem.}
Let $(V,h,Q)$ be a polarized weight two
Hodge structure with $\dim V^{2,0}=1$. Let $K\subset End_{Hod}(V)$
be an imaginary quadratic field such that $V$ is $K$-vector space
and assume $Q(xv,xw)=x\bar{x}Q(v,w)$ for $x\in K$, $v,\,w\in V$.
Let $n=2m=\dim V$.

Then there is a Hodge structure $S$ such that
$$
C^+(Q)\cong S^{2^{m-2}},\qquad\dim_\QQ S=2^{m+1}.
$$
The Hodge structure $S$ splits as:
$$
S\cong S_0\oplus S_1\oplus\ldots\oplus S_m,\qquad S_i\cong S_{m-i},\qquad
\dim_\QQ S_i=2{m\choose i}.
$$
The $S_i$ are simple Hodge structures
except $S_l$ if $2l=m$, $l\equiv 2\;(4)$ and the polarization satisfies an additional condition,
in that case $S_l\cong (S_l')^2$  and $S'_l$ is irreducible.

The $S_i$ are $K$-vector spaces, $K\subset End_{Hod}(S_i)$,
they are simple if $h$ is general and 
$V\hookrightarrow S_0\otimes S_1\;(\subset C^+(Q)\otimes C^+(Q))$. 

\ts To appear.
\qed

\subsection{Example.}(With thanks to M.\ Nori.)\label{nori}
Smooth quartic surfaces in $\PP^3$ are K3 surfaces. Consider a surface
$X$ in $\PP^3$ defined by:
$$
X:\qquad T^4=F(X,Y,Z)
$$
where $F=0$ defines a smooth plane quartic curve $C\;(\subset \PP^2)$.
Let $p\in C$ and let $L_p\;(\subset\PP^2)$ be the tangent line to $C$ at
$L_p$. Then $L_p\cap C=\{p,p',p''\}$ (since $p$ has multiplicty 2).
If $p\in C$ is general and
$u$ is a suitable parameter along $L_p$,
then restriction of $F$ to $L_p$ is given by $u^2(u-1)(u+1)$
(with $u=0$ corresponding to $p\in L_p$).
The curve in $X$ lying over $L_p$ is defined by $T^4=u^2(u^2-1)$,
thus it is irreducible and its
normalization $E_p$
has a $4:1$ map onto $L_p\cong\PP^1$ which is totally ramified
over $u=\pm 1$ and over $u=0$ it has two ramification points. 
Hence $E_p$ is an elliptic curve and since it has an isomorphism of order 4,
$T\mapsto iT,\;i^2=-1$, with fixed points (the points over $u=\pm 1$) one has
$E_p\cong E_i:=\CC/(\ZZ+i\ZZ)$. The elliptic surface
$$
{\cal E}\longrightarrow C,\qquad
{\cal E}:=\cup_{p\in C} E_p\times\{p\}\;\subset\PP^3\times C,
$$
is thus isotrivial and projection on the first factor gives a surjective map
${\cal E}\rightarrow X$. After a suitable base change $C'\rightarrow C$
and normalization the surface will be a product $E_i \times C'$. Thus we have
a surjective rational map:
$$
\pi:E_i\times C'\longrightarrow X,
$$
so $V\;(\subset H^2(X,\QQ))$ lies in the image of $\pi_*$.
Now 
$S_0\cong H^1(E_i,\QQ)$, hence $E_i$ is isogeneous to a simple factor 
$A_0$ of the Kuga-Satake variety $A$ of $V\subset H^2(X,\QQ)$ 
and one verifies that
$S_1\subset H^1(C',\QQ)$, 
hence $Jac(C')$ maps onto a simple factor $A_1$ of $A$.
Then we may take $Z=E_i\times C'$ (or a blow up if $\pi$ is not a morphism), 
the map $\phi$ is the composition:
$$
\phi: E_i\times C'\longrightarrow E'\times Jac(C')
\longrightarrow A_0\times A_1\hookrightarrow A\times A.
$$

\subsection{} 
The following result was inspired by \cite{Vks} (and was generalized following 
a discussion with M. Nori).
It has interesting geometrical applications and
gives a better understanding of the Hodge structures $S_i$.

\subsection{Theorem.} 
Let $(V,h,\Psi)$ be a polarized Hodge structure of weight $k$,
let $K\subset End_{Hod}(V)$
be an imaginary quadratic field such that $V$ is a $K$-vector space
and assume $\Psi(xv,xw)=x\bar{x}\Psi(v,w)$ for $x\in K$, $v,\,w\in V$.

If $K$ acts via scalar multiplication on $V^{k,0}$, then there exist  $(V,h',\Psi')$, $(K,h_K,E_K)$, 
polarized  Hodge structures of weight $k-1$  and one respectively,
such that 
$$
  (V,h)\hookrightarrow (V,h')\otimes (K,h_K),\qquad
  {\rm and}\qquad
  \Psi=(\Psi'\otimes E_K)_{|V}.
$$

\ts A generalization, with $K$ a CM-field, is to appear. \qed


\begin{thebibliography}{DMOS}

\bibitem[CKT]{CKT} J.\ A.\ Carlson, A.\ Kasparian, D.\ Toledo,  
{\it Variations of Hodge structure of maximal dimension},
Duke Math.\ J.\ {\bf 58} 669--694  (1989). 

\bibitem[C]{Ca} J.W.S.\ Cassels, {\it Rational Quadratic Forms},
LMS Monographs {\bf 13}, Academic Press, London, (1978).

\bibitem[D]{DK3} P.\ Deligne, {\it La conjecture de Weil pour les surfaces K3},
Invent.\ Math.\ {\bf 15}  206--226 (1972).

\bibitem[DMOS]{DMOS} P.\ Deligne et al, {\it Hodge Cycles, Motives and Shimura Varieties},  LNM 900, 
Springer-Verlag, Berlin etc.\ (1982).

\bibitem[FH]{FH} W.\ Fulton and J.\ Harris, {\it Representation Theory},
  GTM 129, Springer-Verlag, Berlin etc.\ (1991).

\bibitem[vG]{vG} B.\ van Geemen, {\it An introduction to the Hodge conjecture 
for abelian varieties}. In: Algebraic Cycles and Hodge Theory, LNM 1594,
Springer-Verlag, Berlin etc.\ (1994).

\bibitem[G]{BG} B.\ B.\ Gordon, {\it A survey of the Hodge conjecture for Abelian Varieties}, Duke preprint alg-geom 9709030, to appear in the second edition of `A survey of the Hodge conjecture' by James D.\ Lewis.

\bibitem[GH]{GH} P.\ Griffiths and J.\ Harris,
  {\it Principles of Algebraic Geometry},
John Wiley and Sons, New York (1978).
  

\bibitem[KS]{KS} M.\ Kuga and I.\ Satake, {\it Abelian varieties associated
to polarized $K_3$ surfaces}, Math.\ Ann.\ {\bf 169} 239--242 (1967).

\bibitem[La]{Lam} T.\ Y.\ Lam, {\it The algebraic theory of quadratic forms},
W.\ A.\ Benjamin, Inc., Reading (1973).

\bibitem[Lo]{Lom} G.\ Lombardo, to appear.

\bibitem[M]{M} K.\ Matsumoto, {\it Theta functions on the bounded
symmetric domain of type $I_{2,2}$ and the period map of a 4-parameter
family of K3 surfaces},  
Math.\ Ann.\ {\bf 295}, 383--409 (1993).

\bibitem[Ma]{Ma} R.\ Mayer, {\it Coupled contact systems and rigidity of maximal dimensional variations of
Hodge structure},  Trans. Amer. Math. Soc., to appear (eprint alg-geom/9712001).


\bibitem[MZ1]{MZ} B.J.J.\ Moonen and Yu.G.\ Zarhin, {\it Hodge and Tate classes
on simple abelian fourfolds},
Duke Math.\ J.\ {\bf 77} 553--581 (1995).

\bibitem[MZ2]{MZ2}  B.J.J.\ Moonen and Yu.G.\ Zarhin, {\it Hodge classes on abelian varieties of low dimension}, 
 eprint math/9901113 

\bibitem[P]{Pa} K.\ Paranjape, {\it Abelian varieties associated to certain
    K3 surfaces},
Compos.\ Math.\ {\bf 68} 11--22 (1988). 

\bibitem[S1]{Sata} I.\ Satake, {\it Clifford algebras and families of abelian 
varieities}, Nagoya J.\ Math.\ {\bf 27} 435--446 (1966).

\bibitem[S2]{Sch} W.\ Scharlau, {\it Quadratic and Hermitian Forms},
Grundlehren der Math.\ Wiss.\ 270,  Springer-Verlag, Berlin etc.\ (1985).

\bibitem[S3]{Sab} C.\ Schoen, {\it Hodge classes on self-products of a variety with an automorphism}, Compos.\ Math.\ {\bf 65} 3--32 (1988).

\bibitem[S4]{Scov} C.\ Schoen, {\it Varieties dominated by product varieties},
Int.\ J.\ Math.\ {\bf 7} 541--571 (1996).

\bibitem[S5]{Sdec} C.\ Schoen, {\it Addendum to: hodge classes on self-products of a variety with an automorphism}, Compos.\ Math.\ {\bf 114} 329--336 (1998).

\bibitem[S5]{serre} J.-P.\ Serre, {\it A Course in Arithmetic},
GTM 7, Springer Verlag, Berlin etc.\ (1973).

\bibitem[V]{Vks} C.\ Voisin, 
{\it Remarks on zero-cycles of self-products of varieties.}
M.\ Maruyama (ed.), Moduli of vector bundles. 
(35th Taniguchi
symposium, 1994). Lect. Notes Pure Appl. Math. {\bf 179}, 265--285. 
Marcel Dekker, New York  (1996).

\bibitem[W]{WeilP} A.\ Weil, {\it Variet\`es K\"ahleriennes}.
Hermann, Paris (1971).
\end{thebibliography}
\end{document}